\documentclass[a4paper,12pt]{article}
\font\header=cmssdc10 at 20pt
\usepackage[ansinew]{inputenc}   
\usepackage[brazil]{babel}
\usepackage{amsfonts}
\usepackage{indentfirst}
  \usepackage[normalem]{ulem}
  \usepackage{graphicx}
  \usepackage{amssymb}
 \usepackage{draftcopy}
  \usepackage{amsmath}
  \usepackage{rotating}
  \usepackage{color}
  \usepackage{setspace}
  \usepackage{fancybox}
 \usepackage{fancyvrb}
	\usepackage{fancyhdr}
	\usepackage[top=0.5cm,left=1cm,right=1cm,bottom=0.5cm]{geometry}
  \usepackage{wrapfig}

	\newcommand{\comenta}[1]{%
	} 

  \usepackage{float}
  \usepackage{type1cm}
\usepackage{eso-pic}
\usepackage{color}
 \usepackage{amsthm,amsmath,amsfonts}
\usepackage{fullpage,enumerate}
\usepackage{tikz}
\usepackage[colorlinks,citecolor=blue,urlcolor=blue]{hyperref}
\usepackage{hypernat}
\usepackage[notref,notcite]{showkeys}
\usepackage{graphicx}

\begin{document}

{\header Survival in locally and globally changing environments}

\vskip2cm

Rinaldo B. Schinazi

University of Colorado at Colorado Springs

rschinaz@uccs.edu

\vskip1cm

{\bf Abstract.} We consider branching like models in local, global and fixed environments. We show that survival is more likely in a locally changing environment than in a fixed environment and that survival in a fixed environment is itself more likely than in a globally changing environment.

\vskip1cm

{\header  1 Globally changing environment}

\vskip1cm

 Consider a sequence $(E_i)_{i\geq 1}$ of independent and identical random environments.  Let $(Z_n)_{n\geq 0}$ be a discrete time stochastic process  defined as follows.
The process starts with a single individual. That is,  $Z_0=1$. We sample an  environment $E_1$. This environment determines the offspring distribution $(p_k)_{k\geq 0}$ for this individual. Then we sample an offspring number $k\geq 0$ according to that distribution and let $Z_1=k$.  If $k=0$ then $Z_n=0$ for all $n\geq 1$. If $k\geq 1$ then we sample a new environment  $E_2$ which determines a new offspring distribution. Each of the $k$ individuals in generation 1 samples independently an offspring number from this same offspring distribution. The sum of these numbers is $Z_2$ and so on.  More formally, let $Z_0=1$ and if for $n\geq 1$ let
$$Z_n=\sum_{i=1}^{Z_{n-1}} Y_{n,i}$$
where for fixed $n\geq 1$ the variables $Y_{n,i}$, $i=1,2,\dots$, are independent and identically distributed with distribution $E_i$. For $n\not =m$ and all $i$, $j$, we have that $Y_{n,i}$ and $Y_{m,j}$ are independent but not necessarily identically distributed.

Let $M$ be the mean offspring. Let $P_0$ be the probability of no offspring. The distributions of the random variables $M$ and $P_0$ are the same for all environments $E_i$, $i\geq 1$. Conditional on the offspring distribution $(p_k)_{k\geq 0}$ we have
$$E(M|(p_k))=\sum_{k\geq 0}kp_k$$
and
$$E(P_0|(p_k))=p_0.$$

Assume that 
$E|\ln M|<+\infty$. Smith and Wilkinson (1969) have proved the following.  

$\bullet$ If $E(\ln M)\leq 0$ then the process dies out with probability 1. That is, 
$$P(Z_n\geq 1, \forall n\geq 0)=0.$$

$\bullet$ If $E(\ln M)>0$ and $E|\ln(1-P_0)|<+\infty$ then the process survives. That is,
$$P(Z_n\geq 1, \forall n\geq 0)>0.$$

\medskip

We now give an example.  For each $i\geq 1$ we sample a mean $M_i$ from an uniform distribution in $[0,a]$ where $a>0$ is fixed. The random environment $E_i$ is a Poisson distribution with mean $M_i$. Hence,
$$E|\ln M|=\frac{1}{a}\int_0^a |\ln m| dm<+\infty.$$
We also have
$$E|\ln(1-P_0)|=\frac{1}{a}\int_0^a \ln |1-e^{-m}|dm<+\infty.$$
Therefore, the process survives if and only if
$$E(\ln M)=\frac{1}{a}\int_0^a \ln m dm>0.$$
Computing the integral yields that survival is possible if and only if $\ln a>1$. That is, $a>e$.
On the other hand if we consider a classical branching process (i.e. all offspring distributions $Y_{n,i}$ are i.i.d.) with mean offspring $E(M)=\frac{a}{2}$ then survival is possible for any $a>2$.

\medskip

Going back to the general case we observe that if $M$ is not a constant we have by Jensen's inequality
$$E(\ln M)<\ln E(M).$$
Hence, if $E(\ln M)>0$ then $E(M)>1$. In particular, if the branching process in random environment survives so does a classical branching process with mean offspring $E(M)$. In this sense, the global changing environment hampers survival. 

Haccou and Iwasa (1996) and Haccou and Vatutin (2003) have also compared  global changing and fixed environments for different ecological and biological questions.

\vskip1cm

{\header 2 Locally changing environment}

\vskip1cm

Consider the following discrete time branching process in changing local environments. As in the preceding model we have a sequence $(E_i)_{i\geq 1}$ of independent and identical random environments. Let $r$ be in $[0,1]$. Start the process with a single individual at
time 0. Say that this individual has type 1.
The random environment $E_1$ gives the offspring distribution for type 1. All type 1 individuals give birth independently and with this same distribution. Every time there is a birth there are two possibilities:

$\bullet$  With probability $1-r$, the new individual  has the same type as its parent 
and hence keeps the same offspring distribution as its parent.

$\bullet$ With probability $r$ the new individual is given a new type and hence a new offspring distribution.

We number the types as they appear. 

As for the preceding model let $M$ be the offspring mean. We assume $M$ has a continuous distribution $\mu$ with support contained in $[0,\infty)$ and let $r\in[0,1]$.   
Since $\mu$ is continuous every type appears only once. The population is said to survive if there is a strictly positive probability that at all times there is at least one individual in the population. 
  
 The next result shows that if survival is possible for some $r>0$ then it is possible for all $r<r_c$ where $r_c>0$.

\medskip

{\bf Theorem 1.} 

(i) {\sl If $\mu (m:m>1)>0$ then there exists $r_c>0$ such that the population survives for all $0<r<r_c$.}

(ii) {\sl If $\mu (m:m>1)=0$ the population dies out.} 

\medskip

In particular, the population survives for some $r>0$ if and only if $P(M>1)>0$. 

\medskip

We now give an example. Assume $M$ is uniformly distributed on $[0,a]$ for some fixed $a$.
By Theorem 1 survival is possible in a locally changing environment if and only if $a>1$. 
On the other hand in a fixed environment with mean offspring $E(M)=a/2$, survival is possible 
if and only if $a>2$.

\medskip

Here are two consequences of Theorem 1. Let $M$ be a random variable with distribution $\mu$.

$\bullet$ The model in local random environment can survive for arbitrarily low $E(M)$. This should be compared to the model in global random environment which survives only for $E(\ln M)>0$ and the model in fixed environment which survives only for $E(M)>1$. It is easier to survive in local random environment than in the two other environments. 

$\bullet$ If $\mu$ has infinite support then for any $r<1$ we have $\mu(m:m(1-r)>1)>0$.  Hence, a fixed type has a positive probability of surviving by the argument in the proof of Theorem 1. Therefore, the population survives for all $r$ in $[0,1)$. 

\medskip

We conjecture that there exists a threshold $r_c$ in $(0,1]$ such that survival is possible for all $r<r_c$ and not possible for any $r\geq r_c$. The conjecture for an uniform distribution $\mu$ can be verified by direct computation, see the example below. The conjecture for general $\mu$ is open.

The model with locally changing environment is a discrete time generalization of a model introduced by Cox and Schinazi (2012). The results in this paper show that the phenomena observed
in Cox and Schinazi (2012) for specific offspring distributions are in fact quite general.

\vskip1cm

{\header 3 The tree of types}

\vskip1cm

 The key to our analysis is the so-called  tree of types first introduced in Schinazi and Schweinsberg (2008) for a different model. This tree keeps track of the genealogy of the different types.  Each vertex in the tree will be labeled by a positive integer.  There will be a vertex labeled $k$ if and only if a type $k$ individual is born at some time.  We draw a directed edge from $j$ to $k$ if the first type $k$ individual to be born had a type $j$ individual as its parent.  This construction gives a tree whose root is labeled $1$ because all types are descended from the pathogen of type $1$ that is present at time zero. 

\medskip

{\bf Lemma 1. } {\sl The population survives if and only if the tree of types is infinite.}

\medskip

{\it Proof of Lemma 1}

Note first that if the tree of types is infinite then there are infinitely many births in the population and hence the population survives forever. 

For the converse, assume that the tree of types is finite. Then, there are only finitely many types appearing.  We claim that each type eventually dies out. This is so because if type 1 say does not die out then infinitely many types appear as offspring of type 1 individuals (there are infinitely many births from type 1 parents and each birth has a positive probability $r$ of being of a different type). Hence, if the tree of types is finite, only finitely many types appear and each type dies out. Therefore,  the population dies out. This completes the proof of Lemma 1.

\medskip

We claim that the tree of types is a Galton-Watson tree. This is so because the offsprings of different vertices in the tree of types are independent and identically distributed. Hence, the tree of types is infinite if and only if the mean offspring at a vertex is strictly larger than 1. The next lemma gives a way to compute this mean offspring.

\medskip

{\bf Lemma 2. }{\sl Let $X_1$ be the total number of type 1 individuals that are ever born in the population. Let $X$ be the total number of individuals whose type is not 1 but whose parent is type 1. We have}
$$E(X)=\frac{r}{1-r}E(X_1).$$

\medskip

{\it Proof of Lemma 2}

 It is helpful to represent the population evolution as a tree. The population starts with a single individual represented by the root of the tree. We draw edges between a parent and its children. In order for the population to survive this tree needs to be infinite.

Consider the population starting with a single type 1 individual.  Consider now the subtree including only type 1 individuals and individuals whose parent is type 1. Let $T$ be the total number (possibly infinite) of individuals in this subtree excluding the root. Note that
$$T=X+X_1.$$
Moreover, an individual in the subtree is type 1 with probability $1-r$ and not type 1 with probability $r$. Hence,
$$E(X_1)=(1-r)E(T)\qquad E(X)=rE(T).$$
Therefore,
$$E(X)=\frac{r}{1-r}E(X_1).$$
This completes the proof of Lemma 2.

\medskip

Note that $X$ is the total number of new types that type 1 individuals gave birth to. In particular, $E(X)$ is the mean offspring per vertex in  the tree of types.

\medskip

{\bf Lemma 3. }{\sl Let $E(X)$ be the mean offspring per vertex in the tree of types. A fixed type has a positive probability of surviving forever if and only if
$$E(X)=+\infty.$$}

\medskip

{\it Proof of Lemma 3}

If a fixed type survives forever then it will give birth to infinitely many types (each birth has a fixed probability $r$ of being of a new type). Hence, on the event
of survival we have $X=+\infty$. Therefore,  $E(X)=+\infty$. 
Conversely if $E(X)=+\infty$ then by Lemma 2, $E(X_1)$ is infinite as well. Since the type 1 sub-tree is a Galton-Watson tree it has a positive probability of being infinite.
This completes the proof of Lemma 3.

\medskip

The tree of types analysis shows that we have three distinct behaviors for this model. Let $E(X)$ be the mean offspring in the tree of types.

$\bullet$ If $E(X)\leq 1$ then the population dies out.

\smallskip

$\bullet$ If $1<E(X)<+\infty$ then every type that appears eventually dies out but the population as a whole has a positive probability of surviving.

\smallskip

$\bullet$ If $E(X)=+\infty$ then a fixed type has a positive probability of surviving.

\medskip

We will show below on one example that these three behaviors actually happen.

\vskip1cm

{\header 4 Proof of Theorem 1}

\vskip1cm

{\it Proof of (i)}

Assume  that $\mu (m:m>1)>0.$ Then, there exists $r_c>0$ such that for all $0<r<r_c$ we have
$$\mu(m:m(1-r)>1)>0.$$
Hence, with positive probability type 1 has an offspring distribution with mean $m$ such that $m(1-r)>1$. Consider now only the type 1 individuals in the population. This can be seen as a subtree of the tree representing the whole population. Type 1 individuals all have the same offspring distribution (and hence the same mean offspring $m$). Moreover, the child of a type 1 individual is type 1 with probability $1-r$. Therefore, the mean type 1 offspring of a type 1 individual is $m(1-r)$. This shows that the type 1 subtree is a Galton-Watson process with mean offspring $m(1-r)$. Since $m(1-r)>1$ this sub-tree has a positive probability of being infinite. Hence, the whole population has a positive probability of surviving. This completes the proof of (i).

\medskip

{\it Proof of (ii)}

Assume that $\mu (m:m>1)=0$. Let $m_1\leq 1$ be the mean offspring for type 1 individuals.
Given $m_1$ the subtree of type 1 individuals is a Galton-Watson process with mean offspring $m_1(1-r)\leq 1-r$. Recall that starting with a single individual the offspring of this Galton-Watson process in the first generation is $m_1(1-r)$, in the second generation is $\Big (m_1(1-r)\Big)^2$, in the $k$-th generation is $\Big (m_1(1-r)\Big)^k$ and so on, see for instance Schinazi (2014). We sum over  all generations to get
$$E(X_1|m_1)=\sum_{k=1} \Big (m_1(1-r)\Big)^k=\frac{m_1(1-r)}{1-m_1(1-r)},$$
and by Lemma 2
$$E(X)=rE[\frac{M}{1-M(1-r)}],$$
where $M$ has distribution $\mu$.
Since $M\leq 1$ a.s. we have
$$\frac{M}{1-M(1-r)}\leq \frac{1}{1-(1-r)}=\frac{1}{r}.$$
Hence,
$E(X)\leq r\frac{1}{r}=1.$ That is, the tree of types is finite a.s. and by Lemma 1 the population dies out.
This completes the proof of (ii).

\vskip1cm

{\header 5 An example}

\vskip1cm

Assume that the distribution $\mu$ of $M$ is uniform on $[0,a]$ where $a>0$ is fixed.

$\bullet$  If $a\leq 1$ then $P(M>1)=0$ and the population dies out by Theorem 1.

\smallskip

$\bullet$  If $a(1-r)>1$ then a fixed type has a positive probability of surviving. By Lemma 3, $E(X)=+\infty$ and
the population survives with positive probability.

\smallskip

$\bullet$ If $1<a<2$ and $a(1-r)< 1$ then no fixed type can survive forever. We have
$$E(X)=rE[\frac{M}{1-M(1-r)}]\leq \frac{ar}{1-a(1-r)}<+\infty,$$
and
$$E(X)=\frac{1}{a}\int_0^a \frac{m}{1-m(1-r)}dm.$$
A direct computation shows the existence of $r_c$ in $(1-1/a,1)$ such that if $r<r_c$ then $E(X)>1$ and if $r\geq r_c$
then $E(X)\leq 1$. See Corollary 1 in Cox and Schinazi (2012).

\vskip1cm

{\header References}

\vskip1cm

J.T. Cox and R.B. Schinazi (2012) A Branching Process for Virus Survival. Journal of Applied Probability  49 (2012)  888-894.

P. Haccou and Y. Iwasa (1996) Establishment probabilities in fluctuating environments: a branching model. Theoretical Population Biology 50, 254-280.

P. Haccou and V. Vatutin (2003) Establishment success and extinction in autocorrelated environment. Theoretical Population Biology 64, 303-314.

R.B. Schinazi (2014) {\it Classical and spatial stochastic processes}, second edition. Birkhauser.

R.B. Schinazi and J. Schweinsberg (2008) Spatial and non spatial stochastic models for immune response. Markov Processes and Related Fields14, 255-276.

\end{document}